\documentclass[12pt]{article}
\usepackage{stackrel}
\usepackage{amsmath}
\usepackage{amsfonts}
\usepackage{amssymb}
\usepackage{amsthm}
\usepackage{color}

\newtheorem{thm}{Theorem}[section]

\newtheorem{re}{Remark}[section]

\newcommand{\R}{{\rm I}\kern-0.18em{\rm R}}
\newcommand{\1}{{\rm 1}\kern-0.25em{\rm I}}
\newcommand{\E}{{\rm I}\kern-0.18em{\rm E}}
\newcommand{\p}{{\rm I}\kern-0.18em{\rm P}}


\usepackage{epsfig}

\usepackage{fancyhdr}

\usepackage{graphics}

\usepackage{amsopn}  

\usepackage{amssymb}
\usepackage{amsxtra}
\usepackage{amsfonts}

\usepackage{afterpage}
\usepackage{eucal}
\usepackage{eufrak}

\usepackage{enumerate}

\allowbreak

\def\fnote#1{\footnote}

\newcommand{\bea}{\begin{eqnarray}}

\newcommand{\eea}{\end{eqnarray}}

\newcommand{\beas}{\begin{eqnarray*}}
\newcommand{\eeas}{\end{eqnarray*}}

\title{An autoregressive model leading to stable distributions}
\author{Lev B. Klebanov\footnote{Department of Probability and Statistics,
Charles University, Prague, Czech Republic}, Gregory Temnov\footnote{Department of Probability and Statistics, Charles University, Prague, Czech Republic}, Ashot Kakosyan\footnote{Yerevan State University, Yerevan, Armenia}}
\date{}
\begin{document}
\maketitle

\begin{abstract}
We construct an autoregressive model with random coefficients that has a stationary distribution after proper normalization. This limit distribution is found to be stable. 
\end{abstract}

\section{Introduction}
\setcounter{equation}{0}

Models involving heavy-tailed distributions have been in the focus of researchers in probability and its applications for quite a long time. Unfortunately, only few authors seemed explain convincingly, the reason why heavy-tailed random variables are often observed in many appications. Typical citations of the central limit theorem for convergence to Pareto-L\'{e}vy distributions are not appropriate, as long as this type of convergence is applicable to the summands with heavy tails only. However, several models are knowm, that provide some possible explanations. To cite a few, it is worth to mention some examples from Zolotarev \cite{ZoS}, LePage series representation (see, for example, \cite{ST}) and models by Klebanov, Melamed and Rachev  \cite{KMR}, and by Klebanov and Sl\'{a}mov\'{a} \cite{KS}.  Here we propose an autoregressive model that leads to stable distributions after the proper normalization procedure. Its modification gives another heavy-tailed limit distribution.

\section{First Model}
\setcounter{equation}{0}
Suppose that $G_1,G_2, \ldots , G_n, \ldots$ is a sequence of arrival times of a Poisson process with intensity $1$, and $a>1/2$ is a real number. Define a discrete time random process $\{ X_n, \; n\geq 0 \}$ as follows:
\begin{equation}\label{eq1}
X_o=1, \; X_n = \Bigl(  \frac{G_n}{G_{n-1}}\Bigr)^a X_{n-1} + \varepsilon_n, \; n=1,2, \ldots, 
 \end{equation}
 where the sequence of independent identically distributed (i.i.d.) random variables $\{\varepsilon_n, \; n>0 \}$ is independent also of the Poisson process introduced above. Clearly, the relation (\ref{eq1}) can be considered as an autoregressive model with random coefficients $(G_n/G_{n-1})^a$. From (\ref{eq1}) it follows that
 \[ X_n =\Bigl(  \frac{G_n}{G_{n-1}}\Bigr)^a X_{n-1} + \varepsilon_n=\Bigl(  \frac{G_{n-1}}{G_{n-2}}\Bigr)^a X_{n-2} + \Bigl(  \frac{G_n}{G_{n-1}}\Bigr)^a \varepsilon_{n-1}+\varepsilon_n =\]
 \[ = \ldots = \sum_{k=1}^{n-1}\Bigl(  \frac{G_{n}}{G_{k}}\Bigr)^a \varepsilon_{k}+\varepsilon_n+G_n^a X_o .\]
 After normalizing by $n^a$ we obtain
 \begin{equation}\label{eq2}
 \frac{X_n}{n^a}=\Bigl(\frac{G_n}{n}\Bigr)^a \Bigl( 1+ \sum_{k=1}^{n}\frac{\varepsilon_k}{G_k^a} \Bigr). 
\end{equation}
Convergence to a limit as $n \to \infty$ allows to observe that:
\begin{enumerate}
\item $\Bigl(\frac{G_n}{n}\Bigr)^a \to 1$. Indeed, $G_n=E_1+ \ldots +E_n$, where $E_j$, $j =1, \ldots ,n$ are i.i.d. exponential random variables with unit scale, so that applying the law of large numbers is sufficient.
\item  $\sum_{k=1}^{n}\varepsilon_k /G_k^a \to  \sum_{k=1}^{\infty}\varepsilon_k /G_k^a=Z$, where  $Z$ is the sum of LePage series. If the series converges, this sum has strictly stable distribution (see, for instance, \cite{ST}).
\item From two previous points, we see that the limit distribution of $X_n/n^a$ is a stable distribution.
\end{enumerate}

Next, adding some conditions for the convergence of LePage series, we formulate the following result.
\begin{thm}\label{th1}Let the discrete time process $X_n$ be defined by the relation (\ref{eq1}), where $\{G_n, \; n\geq 1 \}$, $a>1/2$ and $\varepsilon_n$, $n=1, 2, \ldots$ are as mentioned above. Suppose that random variables $\varepsilon_n$ have symmetric distribution, and the moment $\E |\varepsilon_n|^{\alpha}$ exists, $\alpha =1/a$. Then there exists a limit distribution for $X_n/n^a$. This limit distribution is stable with unit location parameter, index of stability $\alpha = 1/a$, skewness $\beta=0$, and scale parameter $\sigma =(\E|\varepsilon_n|^\alpha/c_{\alpha})^a$, where $c_{\alpha}=(1-\alpha)/(\Gamma(2-\alpha)\cos(\pi \alpha/2))$ for $\alpha \neq 1$ and $c_1=2/\pi$.
\end{thm}
\begin{re}\label{re1}
Theorem \ref{th1} gives the sufficient condition for the convergence. However, if we have the convergence of a sequence $\{\varepsilon_n, \; n=1,2,\ldots \}$ then the limit distribution has to be stable (with certain  parameters), which follows from the stability of the sum in LePage series.
\end{re}

\section{Interpretations of the Model}
\setcounter{equation}{0}

1. Suppose that we have a Poisson process with unit intensity; its arrival time are denoted by $G_n$, $(n=1,2, \ldots)$. Consider an electric charge, $X_n$, located at the point $G_n$, with potential defined according to the following rule. The charge at the origin is $X_o=1$; if the charge at point $G_k$ is $X_k$ then the charge at $G_{k+1}$ is
\[ X_{k+1}= \Bigl(  \frac{G_{k+1}}{G_{k}}\Bigr)^2 X_{k} + \varepsilon_k. \]
This rule means that Coulomb's force acting on the charge $X_o$ from $X_{k+1}$ would be  the same in its absolute value as that between $X_o$ and $X_k$ for the case of the absence of the errors $\varepsilon_n$. It is clear that, the farther the charge is, the greater it must be in its absolute value. In this contexy, Theorem \ref{th1} shows that the charge increases at rate $n^2 Y$, where $Y$ is a stable random variable. Its parameters are given in Theorem \ref{th1} for the case $a=2$, that is for $\alpha=1/2$.

2. Similar interpretation can be obtained replacing electrical charges by masses. However, Theorem \ref{th1} does not work in this case. Nevertheless, it is not difficult to see that the convergence is still valid. Therefore, Remark \ref{re1} shows that we have the growth of the mass as $n^2 Z$, where $Z$ is a random variable with L\'{e}vy distribution (with unit location parameter).

3. This model can also have interpretation in term of a collective risk model in which claims arrival times are usually modeled by a Poisson process.  According to the discussed autoregressive model with random coefficients, the surplus value, $X_n$, at a future time moment can be viewed as a future value of the present surplus $X_{n-1}$, summed up with the corresponding claim size, $\varepsilon_n$, due to arrive at time $G_n$.

\end{document}